\documentclass[11pt,reqno]{amsart}

\usepackage{amssymb, amsmath, amsthm, mathrsfs}
\usepackage{hyperref}
\usepackage{enumitem}
\usepackage{float}

\newtheorem*{theoA}{Theorem A}
\newtheorem*{theoB}{Theorem B}
\newtheorem*{theoC}{Theorem C}

\newtheorem*{exm A}{Example A}
\newtheorem*{exm B}{Example B}

\newtheorem*{conjA}{Conjecture A}
\newtheorem*{conjB}{Conjecture B}
\newtheorem*{conC}{Conjecture C}
\newtheorem*{cor A}{Corollary A}
\newtheorem*{cor B}{Corollary B}
\newtheorem{theo}{Theorem}[section]
\newtheorem{lem}{Lemma}[section]

\newtheorem{exm}{Example}[section]

\newtheorem{rem}{Remark}[section]

\newcommand{\be}{\begin{equation}}
\newcommand{\ee}{\end{equation}}
\newcommand{\beas}{\begin{eqnarray*}}
\newcommand{\eeas}{\end{eqnarray*}}
\newcommand{\bea}{\begin{eqnarray}}
\newcommand{\eea}{\end{eqnarray}}

\numberwithin{equation}{section}
\begin{document}
\title [O\MakeLowercase{n the two conjectures}]{O\MakeLowercase{n the two conjectures}}
\date{}
\author[N. S\MakeLowercase{arkar}, , D. P\MakeLowercase{ramanik} \MakeLowercase{and} L. M\MakeLowercase{ahato}]{N\MakeLowercase{abadwip} S\MakeLowercase{arkar}$^*$, D\MakeLowercase{ebabrata} P\MakeLowercase{ramanik} \MakeLowercase{and} L\MakeLowercase{ata} M\MakeLowercase{ahato}}
\address{Department of Mathematics, Raiganj University, Raiganj, West Bengal-733134, India.}
\email{naba.iitbmath@gmail.com}
\address{Department of Mathematics, Raiganj University, Raiganj, West Bengal-733134, India.}
\email{debumath07@gmail.com}
\address{Department of Mathematics, Mahadevananda Mahavidyalaya, Monirampore Barrackpore, West Bengal-700120, India.}
\email{lata27math@gmail.com}

\renewcommand{\thefootnote}{}
\footnote{2020 \emph{Mathematics Subject Classification}: 30D35, 39B32 and 34M10.}
\footnote{\emph{Key words and phrases}: Entire function, sharing values, difference operators.}
\footnote{*\emph{Corresponding Author}: Nabadwip Sarkar.}
\renewcommand{\thefootnote}{\arabic{footnote}}
\setcounter{footnote}{0}

\begin{abstract} In this paper, we investigate the uniqueness problem of entire functions that share an entire function with their higher-order difference operators. We obtain two results that confirm the conjectures posed by Liu and Laine \cite{LL1} and by Zhang et al. \cite{ZKL1}, respectively. In addition, we present several relevant examples to further illustrate and support our findings
\end{abstract}
\thanks{Typeset by \AmS -\LaTeX}
\maketitle
\section{{\bf Introduction and main results}}
In this paper, a meromorphic function $f$ always means a function that is meromorphic in the entire complex plane $\mathbb{C}$. We assume that the reader is familiar with the standard notation and main results of Nevanlinna theory (see, e.g., \cite{WKH1, YY1}).

We define the linear measure by 
\[m(E):=\int_E dt\]
and the logarithmic measure by 
\[l(E):=\int_{E\cap [1,\infty)} \frac{d t}{t}.\]

By $S(r, f)$, we denote any quantity satisfying $S(r, f) = o(T(r, f))$ as $r\to \infty$ possibly outside of an exceptional set of finite logarithmic measure. A meromorphic function $a$ is said to be a \emph{small function} of $f$ if $T(r,a)=S(r,f)$.

Moreover, we use the notations $\rho(f)$, $\rho_1(f)$, $\mu(f)$ and $\lambda(f)$ to denote, respectively, the \emph{order}, \emph{hyper-order}, {emph{lower order}, and \emph{the exponent of convergence of zeros} of a meromorphic function $f$.
If $\mu(f)=\rho(f)$, we say that $f$ is of \emph{regular growth}. 

It is well known that if $f=e^g$, where $g$ is a polynomial, then $\rho(f)=\mu(f)=\deg(g)$. Furthermore, if $f$ and $g$ are two non-constant meromorphic functions in $\mathbb{C}$ such that $\rho(f)<\mu(g)$, then 
\[T(r,f)=o(T(r,g)\;\;(r\to\infty)\]
(see Theorem 1.18 \cite{YY1}).

\smallskip
The \emph{Borel exceptional} small function $a$ of $f$ is defined by
\[\lambda(f-a)=\limsup\limits_{r\to\infty} \frac{\log^+ N(r,0;f-a)}{\log r}<\rho(f)\]
where $\lambda(f-a)$ is the exponent of convergence of zeros of $f(z)-a(z)$. Also we know that if $f$ is a transcendental meromorphic function with finite positive order $\rho(f)$, then $f$ has at most two \emph{Borel exceptional} values (see \cite[Theorem 2.10]{YY1}).

As usual, the abbreviations CM and IM stand for counting multiplicities and ignoring multiplicities, respectively.

Let $c\in \mathbb{C}\backslash \{0\}$. Then the forward difference $\Delta^n_c f$ for each integer $n\in\mathbb{N}$ is defined in the standard way by 
\[\Delta^1_c f(z)=\Delta_c f(z)=f(z+c)-f(z)\]
\[\Delta_c^n f(z)=\Delta_c\left(\Delta_c^{n-1}f(z)\right)=\Delta_c^{n-1} f(z+c)-\Delta_c^{n-1}f(z),\;\;n\geq 2.\]
Moreover
\[ \Delta_c^n f(z)=\sum\limits_{j=0}^n(-1)^{n-j}C^{j}_{n}f(z+jc),\]
where $C^{j}_{n}$ is a combinatorial number. If an equation includes shifts or differences of $f$, then the equation is called difference equation.

\smallskip
We recall the following conjecture proposed by Br\"{u}ck \cite{RB1}.
\begin{conjA} Let $f$ be a non-constant entire function such that $\rho_1(f)\not\in\mathbb{N}\cup\{\infty\}$.
If $f$ and $f'$ share one finite value $a$ CM, then $f'-a=c(f-a),\;\text{where}\;\;c\in\mathbb{C}\backslash \{0\}$. 
\end{conjA}

While the conjecture remains unsettled in its full generality, it has generated a substantial body of research concerning the uniqueness of entire and meromorphic functions sharing a single value with their derivatives.

The study of meromorphic solutions of complex difference equations, as well as the value distribution and uniqueness of complex differences, has become an area of active research. This development is based on the Nevanlinna value distribution theory for difference operators established by Halburd and Korhonen \cite{HK1}, and independently by Chiang and Feng \cite{CF1}. Recently, several authors (see \cite{MBA1}-\cite{BA}, \cite{BM1}, \cite{BM2}, \cite{BM3}, \cite{CC2}-\cite{ZXC1}, \cite{DM1}, \cite{DM2}, \cite{HKLR}-\cite{LLLX1}, \cite{M1}, \cite{MD}, \cite{MD1}, \cite{MP1}, \cite{MS1}, \cite{MS2}, \cite{MNS1}, \cite{MNS}, \cite{QLY1}, \cite{WC1}  \cite{ZKL1}) have investigated the problems concerning value sharing between meromorphic functions and their difference operators or shifts. It is also well known that the operator $\Delta_c f$ can be regarded as the difference analogue of $f'$. We now recall the following result due to Heittokangas et al. \cite{HKLRZ}, which represents a difference analogue of the Br\"{u}ck conjecture.

\begin{theoA}\cite[Theorem 1]{HKLRZ} Let $f$ be a non-constant meromorphic function with $\rho(f)<2$ and $c\in\mathbb{C}$. If $f(z)$ and $f(z+c)$ share the values $a\in\mathbb{C}$ and $\infty$, then $f(z+c)-a=\tau (f(z)-a)$ holds for some constant $\tau$.
\end{theoA}

In the same paper, they also provided the following example to demonstrate that the condition ``$\rho(f)<2$'' cannot be relaxed to ``$\rho(f)\leq 2$''.

\begin{exm A} Let $f(z)=e^{z^2}+1$ and $c\in\mathbb{C}\backslash \{0\}$. Clearly $f(z)$ and $f(z+c)$ share $1$ and $\infty$ CM, 
\[\frac{f(z+c)-1}{f(z)-1}=e^{2cz+c^2}\neq \;\text{constant}.\]
\end{exm A}

For a transcendental entire function of order $\rho(f)<2$, Liu and Laine \cite{LL1} obtained the following result
by replacing $f(z+c)$ by $\Delta^n_cf$ in \emph{Theorem A}.

\begin{theoB}\cite[Theorem 1.3]{LL1} Let $f$ be a transcendental entire function of order $\rho(f)<2$ not
having period $c$. If $f$ and $\Delta^n_cf$ share the value $0$ CM, then $\Delta^n_cf=\tau f$, where $\tau\in\mathbb{C}\backslash \{0\}$.
\end{theoB}

Next we consider the following example.
\begin{exm B}\cite{LL1} Let 
\[f(z)=Ae^{z\log(c+1)}-((1-c)/c),\]
where $c\in\mathbb{R}\backslash \{0\},c>-1$ and $A$ is an arbitrary constant. Then $(\Delta_1f(z)-1)=c(f(z)-1)$.
\end{exm B}

Therefore \emph{Example A} suggests that 
\[\frac{\Delta^nf-a}{f-a}\] 
may reduce to a non-zero constant, at least if $\rho(f)=1$ and $N(r,1/f)\neq S(r,f)$.

Motivated by \emph{Example B}, Liu and Laine \cite{LL1} proposed the following conjecture (see Remark 1.4 \cite{LL1}).

\begin{conjB} Let $f$ be a transcendental entire function of order $1<\rho(f)<\infty$ not
having period $c$. If $f$ and $\Delta^n_cf$ share the value $a\in\mathbb{C}\backslash \{0\}$ CM, then $\Delta^n_cf-a=\tau (f-a)$, where $\tau\in\mathbb{C}\backslash \{0\}$.
\end{conjB}
 
In 2015, Zhang et al. \cite{ZKL1} partially resolved Conjecture B under the assumptions ``$\rho(f)<2$'' and ``$\lambda(f-a)<\rho(f)$''
In fact, they established the following result.

\begin{theoC}\cite[Theorem 1.2]{ZKL1} Let $f$ be a transcendental entire function such that $\rho(f)< 2$ and $\alpha(\not\equiv 0)$ be an entire function such that  $\rho(\alpha)<\rho(f)$ and $\lambda(f-\alpha)<\rho(f)$. If $f-\alpha$ and $\Delta^n_1f-\alpha$ share $0$ CM, then $\alpha$ is a polynomial of degree at most $n-1$ and $f$ must be of form $f(z)=\alpha(z)+H(z)e^{dz}$, where $H$ is a polynomial such that $cH=-\alpha$ and $c, d\in\mathbb{C}\backslash \{0\}$ such that $c^d=1$.
\end{theoC}

In the same paper, Zhang et al. \cite{ZKL1} made the following remark.\par
``We are not sure whether the assumption $\alpha\not\equiv 0$ in Theorem 1.2 \cite{ZKL1} is necessary or not.'' 

Consequently, Zhang et al. \cite{ZKL1} proposed the following conjecture.

\begin{conC} Let $f$ be a transcendental entire function such that $\rho(f)<2$ and $\lambda(f)<\rho(f)$. If $f$ and $\Delta^nf$ share $0$ CM, then $f$ must be of the form $f(z)=ce^{dz}$, where $c,\;d\in\mathbb{C}\backslash \{0\}$.
\end{conC}

To the best of our knowledge, \emph{Conjectures} B and C remain unconfirmed. In this paper, we confirm both conjectures and further improve \emph{Theorems} B and C by eliminating the condition ``$\rho(f)<2$'' 

\subsection{\bf{In connection with Conjecture B}}

\medskip
Concerning \emph{Conjecture B}, we obtain the following result.

\begin{theo}\label{t2} Let $f$ be a transcendental entire function such that $\rho(f)<+\infty$, $c\in\mathbb{C}\backslash \{0\}$ such that $\Delta_c^nf\not\equiv 0$ and $\lambda(f)<\rho(f)$. If $f$ and $\Delta^n_cf$ share $0$ CM, then $f$ must be of form $f(z)=c_0e^{dz}$, where $c_0, d\in\mathbb{C}\backslash \{0\}$.
\end{theo}  

The following example demonstrates that the condition ``$\lambda(f)<\rho(f)$'' in Theorem \ref{t2} cannot be weakened.

\begin{exm} 
Let $f(z)=(\exp z-1)\exp\left(\frac{\log (1+\tau)}{c}z\right)$, where $\log $ denotes the principal branch of the logarithm and $c=2\pi i$ such that $\log (1+\tau)\neq c$. Note that 
\beas \Delta_c f(z)&=&(\exp z-1)\exp\left(\frac{\log (1+\tau)}{c}(z+c)\right)-(\exp z-1)\exp\left(\frac{\log (1+\tau)}{c}z\right)\\
&=&(\exp z-1)\exp\left(\frac{\log (1+\tau )}{c}z\right)\left(\exp(\log (1+\tau ))-1\right)\\
&=&\tau (\exp z-1)\exp\left(\frac{\log (1+\tau )}{c}z\right)\\&=& \tau f(z).\eeas
Clearly $f$ and $\Delta_c f$ share $0$ CM. On the other hand, we see that $\rho(f)\leq 1$ and $\lambda(f)=\lambda (\exp z-1)=1$. Since $\lambda(f)\leq \rho(f)$, it follows that $\lambda(f)=\rho(f)$. Also it is clear that $f$ does not satisfies the conclusion of Theorem \ref{t2}.
\end{exm}

The following example indicates that the condition ``$\rho(f)<+\infty$'' in Theorem \ref{t2} is the best possible.

\begin{exm} 
Let $f(z)=e^z e^{s(z)}$, where $s(z)$ is a periodic function with period $c=\log 2$. Clearly $\rho(f)=+\infty$. Note that $\Delta_c f=f$ and so $f$ and $\Delta_c f$ share $0$ CM. On the other hand, we see that $\lambda(f)=0<\rho(f)$, but $f$ does not satisfies the conclusion of Theorem \ref{t2}.
\end{exm}

\subsection{\bf{In connection with Conjecture C}}

\medskip
Concerning \emph{Conjecture C}, we obtain the following result.

\begin{theo}\label{t1} Let $f$ be a transcendental entire function such that $\rho(f)<+\infty$, $c\in\mathbb{C}\backslash \{0\}$ such that $\Delta_c^nf\not\equiv 0$ and $\alpha(\not\equiv 0)$ be an entire function such that  $\rho(\alpha)<\rho(f)$ and $\lambda(f-\alpha)<\rho(f)$. If $f-\alpha$ and $\Delta^n_cf-\alpha$ share $0$ CM, then
$\alpha$ is a polynomial of degree at most $n-1$ and $f$ must be of form $f(z)=\alpha(z)+H(z)e^{dz}$, where $H$ is a polynomial such that $c_0H=-\alpha$ and $c_0, d\in\mathbb{C}\backslash \{0\}$ such that $c^d=1$.
\end{theo}

The following examples verify the sharpness of the condition `$\lambda(f-\alpha)<\rho(f)$'' in Theorem \ref{t1}.

\begin{exm} Let $f(z)=e^z$ and $c=\log 2$. Note that 
\beas \Delta_c^nf(z)=\sum\limits_{j=0}^n(-1)^jC^j_n f(z+(n-j)c)&=&e^z\sum\limits_{j=0}^n(-1)^jC^j_ne^{(n-j)c}\\&=&
\left(e^{nc}-C^1_ne^{(n-1)c}+\ldots+(-1)^n\right)e^z\\&=&\left(e^{c}-1\right)^n e^z=e^z.
\eeas
Therefore $\Delta_c^nf\equiv f$ and so $f$ and $\Delta_c^nf$ share $1$ CM. Since $0$ and $\infty$ are the Borel exceptional values of $f$, it follows that $1$ can not be a Borel exceptional value of $f$ and so $\lambda(f-1)=\rho(f)$. Clearly $f$ does not satisfies the conclusion of Theorem \ref{t1}.
\end{exm}

\begin{exm} Let $f(z)=\sin z$. Note that $\lambda(f)=\rho(f)$ and $\Delta_{\pi}(f)$ and $f$ share $0$ CM, but $f$ does not satisfy Theorem \ref{t2}.
\end{exm}

It is easy to see that condition ``$\rho(\alpha)<\rho(f)$'' in Theorem \ref{t1} is sharp.

\begin{exm} Let $f(z)=e^z+1$, $\alpha(z)=e^z-1$ and $c=\log 2$. Note that $\rho(\alpha)=\rho(f)$ and $\Delta_cf(z)=e^z$. Clearly $\lambda(f-\alpha)=0<\rho(f)$ and $f-\alpha$ and $\Delta_cf-\alpha$ share $0$, but $f$ does not satisfies the conclusion of Theorem \ref{t1}.
\end{exm}

Following example shows that condition ``$\rho(f)<+\infty$'' in Theorem \ref{t1} is necessary.

\begin{exm} 
Let $f(z)=e^z\left(e^{s(z)}-1\right)$, where $s(z)$ is a periodic function with period $c=\log 2$ and $\alpha(z)=-e^z$. Clearly $\rho(f)=+\infty$. Note that $\Delta_c f=f$ and so $f-\alpha$ and $\Delta_c f-\alpha$ share $0$ CM. On the other hand, we see that $\lambda(f-\alpha)=0<\rho(f)$, but $f$ does not satisfies the conclusion of Theorem \ref{t1}.
\end{exm}

\section {{\bf Auxiliary lemmas}}
For the proof of our main results, we make use of the following five key lemmas.

\begin{lem}\label{l4} \cite[Corollary 2.6]{CF1} Let $f$ be a meromorphic function of finite order $\rho$ and $\eta_1,\;\eta_2\in\mathbb{C}$ such that $\eta_1\neq \eta_2$. Then for any $\varepsilon >0$ we have
\beas m\left(r,\frac{f(z+\eta_1)}{f(z+\eta_2)}\right)=O(r^{\rho-1+\varepsilon }).\eeas
\end{lem}

\begin{lem}\label{l7} \cite[Corollary 8.3]{CF1} Let $\eta_1,\eta_2$ be two arbitrary complex numbers and $f$ be a meromorphic function of finite order $\sigma $. Let $\varepsilon >0$ be given. Then there exists a subset $E\subset (0,\infty)$ with finite logarithmic measure such that for all $r\not\in E\cup[0,1]$ we have
\beas \exp(-r^{\sigma -1+\varepsilon })\leq \bigg|\frac{f(z+\eta_1)}{f(z+\eta_2)}\bigg|\leq \exp(r^{\sigma -1+\varepsilon }).\eeas
\end{lem}

Now we recall the definition of difference polynomial (see \cite[pp. 557]{ILY1}). Let $f$ be a non-constant meromorphic function, $\alpha_{\lambda}\in S(f)$, and let $\mu_{\lambda,j}\geq 0$ and $\tau_{\lambda}\geq 1$ be integers, where $j=1,\ldots, \tau_{\lambda}$ and $\lambda\in J$, a finite index set. We always refer to an expression of the following form when we say ``difference product''
\[M_{\lambda}(z,f)=\prod\limits_{j=1}^{\tau_{\lambda}}\left(f(z+c_{\lambda,j})\right)^{\mu_{\lambda,j}},\]
where $c_{\lambda,1}, \ldots, c_{\lambda,\tau_{\lambda}}$ are finite complex numbers such that at least one of $c_{\lambda,1}, \ldots, c_{\lambda,\tau_{\lambda}}$ is non-zero. Then $M(z,f)=\sum\limits_{\lambda\in J}a_{\lambda} M_{\lambda}(z,f)$ is called a
difference polynomial in $f$. The maximal total degree of $M(z,f)$ in $f$ and the shifts of $f$ is defined by
\[\deg(M(z,f))=\max\limits_{\lambda\in J}\sum\limits_{j=1}^{\tau_{\lambda}} \mu_{\lambda,j}.\]

\begin{lem}\label{l5}\cite[Theorem 2.3]{ILY1}
Let $f$ be a transcendental meromorphic solution of finite order $\rho$ of a difference equation of the form 
$U(z,f)P(z,f)=Q(z,f)$,
where $U(z,f), P(z,f), Q(z,f)$ are difference polynomials such that the total degree $\deg\left(U(z,f)\right)=n$ in $f$ and its shifts and $\deg\left(Q(z,f)\right)\leq n$. Moreover we assume that $U(z,f)$ contains just one term of maximal total degree in $f$ and its shifts. Then for each $\varepsilon>0$ we have 
$$m(r,P(z,f))=O(r^{\rho-1+\varepsilon })+S(r,f)$$
possible outside of an exceptional set of finite logarithmic measure.
\end{lem}

\begin{rem}\label{r1} Also from the proof of Lemma \ref{l5}, we can say that if the coefficients of
$U(z,f)$, $P(z,f)$ and $Q(z,f)$, say $a_{\lambda}$ satisfy $m(r,a_{\lambda})= S(r,f)$, then the same conclusion still holds.
\end{rem}

\begin{lem}\label{l8}\cite[Corollary 3.2]{LLLX1} Let $g$ be a non-constant meromorphic solution of the linear difference equation 
\[b_kg(z+kc)+b_{k-1}g(z+(k-1)c)+\ldots+b_0g(z)=R(z),\]
where $R$ is polynomial and $b_i$'s for $i=0,1,\ldots, k$ are complex constants with $b_kb_0\neq 0$, $c\in \mathbb{C}\backslash \{0\}$ and $k\in \mathbb{N}$. Then either $\rho(g)\geq 1$ or $g$ is a polynomial. In particular if $b_k\neq \pm b_0$, then $\rho(g)\geq 1$.
\end{lem}

\begin{lem}\label{l3}\cite[Theorem 1.51]{YY1} Suppose that $f_1, f_2 ,\ldots, f_n\; (n\geq 2)$ are meromorphic functions and $g_1, g_2, \ldots, g_n$ are entire functions satisfying the following conditions 
\begin{enumerate}
\item[(i)] $\sum\limits_{j=1}^nf_je^{g_j}=0$
\item[(ii)] $g_i-g_j$ is non-constant for $1\leq i<j\leq n$;
\item[(iii)] $T(r,f_j)=o\left(T(r,e^{g_h-g_k})\right)$ $(r\rightarrow \infty,r\not\in E)$ for $1\leq j\leq n$, $1\leq h < k\leq n$.
\end{enumerate}

Then $f_j\equiv 0$ for $j=1,2,\ldots,n$.
\end{lem}

\section {{\bf Proofs of main results}} 
\begin{proof}[{\bf Proof of Theorem \ref{t1}}]
By the given conditions, we have $\lambda(f-\alpha)<\rho(f)$. By Hadamard factorization theorem, there exist an entire function $H(\not\equiv 0)$ and a polynomial $P$ such that 
\bea\label{al.1} f=\alpha+He^{P},\eea
where $\rho(H)<\rho(f-\alpha)$ and $\deg(P)=\rho(f-\alpha)$. Since $\rho(\alpha)<\rho(f)$, we have
$\rho(H)<\rho(f)$ and $\deg(P)=\rho(f)$. Let 
\bea\label{al.1a} P(z)=a_sz^s+a_{s-1}z^{s-1}+\cdots+a_0,\eea
where $a_s(\neq 0),a_{s-1},\ldots,a_0\in\mathbb{C}$ and $s\in\mathbb{N}$. Also from (\ref{al.1}), we get
\bea\label{al.2} \Delta_c^nf=\Delta_c^n\alpha+H_ne^{P},\eea
where 
\bea\label{al.3} H_n(z)=\sum\limits_{i=0}^nc_iH(z+ic)e^{P(z+ic)-P(z)}\;\;\text{and}\;\;c_i=(-1)^{n-i}C^{i}_{n}.\eea

\medskip
First we suppose that $\rho(f)<2$. Then by \emph{Theorem C}, we conclude that $\alpha$ is a polynomial of degree at most $n-1$ and $f$ must be of the form $f(z)=\alpha(z)+H(z)e^{dz}$, where $H$ is a polynomial such that $c_0H=-\alpha$ and $c_0, d\in\mathbb{C}\backslash \{0\}$ such that $c^d=1$.\par

\medskip
Next we suppose that $\rho(f)\geq 2$. Since $f-\alpha$ and $\Delta_c^nf-\alpha$ share $0$ CM, then there exists a polynomial function $Q$ such that 
\bea\label{al.4}\Delta_c^nf-\alpha=(f-\alpha)e^{Q}.\eea

We know that $\rho(H)<\rho(f)$ and so $\rho(H(z+ic))<\rho(f)$ for $i=0,1,\ldots, n$. Note that $\deg\left(P(z+ic)-P(z)\right)\leq s-1=\rho(f)-1$. Then from (\ref{al.3}), we get $\rho(H_n)<\rho(f)$. Also we know that $\rho(\Delta^n_c \alpha)\leq \rho(\alpha)$. Consequently from (\ref{al.2}), we get $\rho\left(\Delta^n_c f\right)\leq \rho(f)$. Now from (\ref{al.4}), we have $\deg(Q)\leq \rho(f)=s$. 

We now consider following two cases.\par

\medskip
{\bf Case 1.} Let $\deg(Q)=0$. Then from (\ref{al.4}), we may assume that
\bea\label{al.5} \Delta_c^nf-\alpha=c_0(f-\alpha),\eea
where $c_0\in\mathbb{C}\backslash \{0\}$. Now from (\ref{al.1}), (\ref{al.2}) and (\ref{al.5}), we get
\bea\label{al.6} \Delta^n_c\alpha-\alpha=(c_0H-H_n)e^P.\eea

Note that $\rho\left(\Delta^n_c\alpha-\alpha\right)<\rho(f)=\rho\left(e^P\right)=\mu\left(e^P\right)$
and $\rho\left(c_0H-H_n\right)<\rho(f)=\mu\left(e^P\right)$. Consequently $T\left(r,\Delta^n_c\alpha-\alpha\right)=S\left(r,e^P\right)$ and $T\left(r,c_0H-H_n\right)=S\left(r,e^P\right)$. Now using Lemma \ref{l3} to (\ref{al.6}), we have
\[\Delta^n_c\alpha\equiv \alpha\;\;\text{and}\;\;c_0H\equiv H_n.\] 

Let $h(z)=e^{sca_sz^{s-1}}$. Then $\rho(h)=\mu(h)$ and
\bea\label{al.7}T(r,h)=m(r,h)=\frac{|s||c||a_s|}{\pi}(1+o(1))r^{s-1}.\eea

Now for $0\leq j\leq n$, we have 
\beas e^{P(z+jc)-P(z)}=e^{jsca_sz^{s-1}}e^{P_j(z)}=h^j(z)e^{P_j(z)},\eeas
where $P_j$ is polynomial with degree at most $s-2$. 
Note that 
\[\rho\left(e^{P_j}\right)=\deg(P_j)\leq s-2<s-1=\rho(h)=\mu(h)\]
and so $T\left(r, e^{P_j}\right)=S(r,h)$ and so $T\left(r, e^{P_j}\right)=S(r,h^j)$ for $j=0,1,2,\ldots, n$.
Now from Lemma \ref{l4}, we have 
\bea\label{al.8} m\left(r,\frac{H(z+jc)}{H(z)}\right)=O\left(r^{\rho (H)-1+\varepsilon }\right),\eea
where $\varepsilon >0$ is arbitrary. Since $\rho(H)<\rho(e^{P})$, we chose $\varepsilon>0$ such that 
\[\rho (H)-1+2\varepsilon < \rho(e^P)-1=\rho(h)=\mu(h).\]

On the other hand, by the definition of the lower order, there exists a $R>0$ such that 
\[T(r,h)>r^{\mu(h)-\varepsilon}=r^{\rho(h)-\varepsilon}\]
holds for $r\geq R$. Then from (\ref{al.8}), we get
\[\lim\limits_{r\to\infty}\frac{m\left(r,\frac{H(z+jc)}{H(z)}\right)}{T(r,h)}=0\]
 and so 
\bea\label{al.9} m\left(r,\frac{H(z+jc)}{H(z)}\right)=S(r,h)\eea
for $j=0,1,\ldots,n$.
Similarly we can prove that
\bea\label{al.10} m\left(r,\frac{H(z)}{H(z+jc)}\right)=S(r,h)\eea
for $j=0,1,\ldots,n$.
Let
\bea\label{al.11} b_{n-j}(z)=c_j\frac{H(z+jc)}{H(z)}e^{P_j(z)},\eea
for $j=0,1,2,\ldots, n$ and  
\bea\label{al.12} F_n(h)=\sum\limits_{j=0}^nb_{n-j}h^j.\eea

Since $T\left(r, e^{P_j}\right)=S(r,h)$, for $j=0,1,\ldots,n$, from (\ref{al.9}), (\ref{al.10}) and (\ref{al.11}), we see that
\bea\label{al.13} m(r,b_j)+m\left(r,\frac{1}{b_j}\right)=S(r,h),\;\text{for}\;j=0,1,\ldots,n.\eea

Now from the proof of Theorem 2.1 \cite{LLLX1}, we deduce that
\beas n\;m(r,h)+S(r,h)=m(r,F_n(h))+S(r,h)\eeas
and so from (\ref{al.7}), we have   
\bea\label{al.18} m(r,F_n(h))=n\frac{|s||c||a_s|}{\pi}(1+o(1))r^{s-1}+S(r,h)=n\;T(r,h)+S(r,h).\eea

Since $c_0H\equiv H_n$, from (\ref{al.3}), we have
\beas F_n(h)=\sum\limits_{j=0}^nb_{n-j}h^j=c_0\eeas
and so $m(r,F_n(h))=O(1)$. Then from (\ref{al.18}), we get
$T(r,h)=S(r,h)$, which is impossible.\par

\medskip
{\bf Case 2.} Let $\deg(Q)\geq 1$. Now by taking logarithmic differentiation on (\ref{al.4}), we get
\beas\frac{(\Delta_c^n f)'-\alpha '}{\Delta f-\alpha}-\frac{f'-\alpha '}{f-\alpha}=Q'\eeas
and so from (\ref{al.1}) and (\ref{al.2}), we obtain
\bea\label{al.19}\left(H_n'H-H_nH'-H_nHQ'\right)e^{P}=(\Delta_c^n\alpha-\alpha)\left(H'+HP'+HQ'\right)-(\Delta\alpha-\alpha)'H.\eea

We consider following two sub-cases.\par

\medskip
{\bf Sub-case 2.1.} Let $\Delta_c^n\alpha\not=\alpha$. We know that $\rho(H)<\rho(e^P)$, $\rho(H_n)<\rho(e^P)$ and $\rho(\alpha)<\rho(e^P)$. Therefore using Lemma \ref{l3} to (\ref{al.19}), we obtain
\bea\label{al.20} H_n'H-H_nH'-H_nHQ'=0\eea
and 
\bea\label{al.21} (\Delta_c^n\alpha-\alpha)(H'+HP'+HQ')-(\Delta_c^n\alpha-\alpha)'H=0.\eea

Then from (\ref{al.20}), we get
\beas \frac{H_n'}{H_n}-\frac{H'}{H}-Q'=0\eeas
and so from (\ref{al.21}), we have 
\beas \frac{(\Delta_c^n\alpha-\alpha)'}{(\Delta_c^n\alpha-\alpha)}=\frac{H_n'}{H_n}+P'.\eeas

On integration, we have
\bea\label{al.22}\Delta_c^n\alpha-\alpha=c_1 H_ne^{P}\;\;(c_1\in\mathbb{C}\backslash \{0\}).\eea

Since $\rho(H_n)<\rho(e^P)$, it follows that $\rho\left(H_ne^{P}\right)=\rho(e^P)$. Also we know that $\rho\left(\Delta_c^n\alpha-\alpha\right)<\rho(e^P)$. Therefore, from (\ref{al.22}), we get a contradiction.\par

\medskip
{\bf Sub-case 2.2.} Let $\Delta_c^n\alpha=\alpha$. Now from (\ref{al.19}), we have 
\[\frac{H_n'}{H_n}-\frac{H'}{H}-Q'\equiv 0\]
and so $H_n\equiv c_2He^Q$, where $c_2\in\mathbb{C}\backslash \{0\}$. Also from (\ref{al.4}), we have
\bea\label{al.23} F_n(h)=\sum\limits_{j=0}^nb_{n-j}h^j=e^{Q},\eea
where $b_{j}$'s are given by (\ref{al.11}). Now from (\ref{al.23}), we conclude that $F_n(h)$ is an entire function and so \[T(r,F_n(h))=m(r,F_n(h)).\]

Then from (\ref{al.18}), we get 
\[T(r,F_n(h))=n\;T(r,h)+S(r,h)\]
and so from (\ref{al.23}), we have 
\[n\;T(r,h)+S(r,h)=T(r,e^Q).\]

Clearly $\rho(h)=\rho(e^Q)$. Since $\rho(h)=s-1$, it follows that $\rho(e^Q)=s-1$ and so $\deg(Q)=s-1$. 

Let 
\bea\label{al.23a} Q(z)=d_{s-1}z^{s-1}+d_{s-2}z^{s-2}+\cdots+d_0.\eea

Now from (\ref{al.4}), we have 
\bea\label{al.24} \sum\limits_{j=1}^n c_j\frac{H(z+jc)}{H(z)}e^{R_j(z)}+(-1)^n-e^{Q(z)}=0,\eea
where 
\[R_j(z)=P(z+jc)-P(z)\;(j=1,2,\ldots,n).\]

Then from (\ref{al.1a}), we may assume that
\bea\label{al.24a} R_j(z)=jsa_scz^{s-1}+P_{s-2,j}(z),\eea
where $P_{s-2,j}(z)$ is a polynomial with degree at most $s-2$. Clearly $\deg(R_j)=s-1$ for $j=1,2,\ldots,n$.

\medskip
First we suppose that $n=1$. Then from (\ref{al.24}), we have
\bea\label{al.25} c_1\frac{H(z+c)}{H(z)}e^{R_1(z)}-1=e^{Q(z)}.\eea

It is clear from (\ref{al.25}) that $\frac{H(z+c)}{H(z)}$ is an entire function. Now from (\ref{al.8}), we have
\[T\left(r,\frac{H(z+c)}{H(z)}\right)=m\left(r,\frac{H(z+c)}{H(z)}\right)=O\left(r^{\rho(H)-1+\varepsilon}\right)\]
and so 
\[\rho\left(\frac{H(z+c)}{H(z)}\right)=\rho(H)-1<\rho(f)-1=s-1=\rho(e^{R_1}).\]
 
Therefore it is easy to verify that $0$ is a Borel exceptional value of the entire function 
\[c_1\frac{H(z+c)}{H(z)}e^{R_1(z)}.\]

Consequently $1$ is not a Borel exceptional of 
\[c_1\frac{H(z+c)}{H(z)}e^{R_1(z)}\] 
and so 
\[c_1\frac{H(z+c)}{H(z)}e^{R_1(z)}-1\]
must have infinitely many zeros. Therefore, from (\ref{al.25}), we get a contradiction.

\medskip
Next we suppose that $n\geq 2$. Then from (\ref{al.23a}) and (\ref{al.24a}), we see that
\bea\label{al.25a} R_j(z)-Q(z)=(jsa_sc-d_{s-1})z^{s-1}+\ldots,\eea
where $j=1,2,\ldots,n$.

We now proceed by considering the following two sub-cases.\par

\medskip
{\bf Sub-case 2.2.1.} Let $j_0sa_sc=d_{s-1}$, where $1\leq j_0\leq n$. Then from (\ref{al.25a}), we can say that $\deg(R_{j_0}-Q)\leq s-2$. In this case, from (\ref{al.24}), we get 
\bea\label{al.26}\left(\sum\limits_{\substack{1\leq j\leq n\\j\neq j_0}}c_j\frac{H(z+cj)}{H(z)}e^{P(z+jc)-P(z+c)}+B_{j_0}e^{P(z+j_0c)-P(z+c)}\right)e^{R_1(z)}=(-1)^{n+1},\eea
where 
\bea\label{al.27} B_{j_0}(z)=c_{j_0}\frac{H(z+j_0c)}{H(z)}-e^{Q(z)-R_{j_0}(z)}.\eea
Let $Q_1=e^{R_1}$. Note that
\beas Q_1(z+(j-1)c)Q_1(z+(j-2)c)\ldots Q_1(z+c)&=&e^{\left(\sum\limits_{i=2}^{j}{P(z+ic)-P(z+(i-1)c)}\right)}\\
&=&e^{P(z+jc)-P(z+c)}\eeas
for $j=2,3,\ldots,n.$ Then (\ref{al.26}) can be written as
\bea\label{al.28} U(z,Q_1)Q_1=(-1)^{n+1},\eea
where 
\beas U(z,Q_1(z))&=&\sum\limits_{\substack{1\le j\leq n\\j\neq j_0}}c_{j}\frac{H(z+jc)}{H(z)}Q_1(z+(j-1)c)Q_1(z+(j-2)c)\cdots Q_1(z+c)\\
  &&+B_{j_0}(z)Q_1(z+(j_0-1)c)Q_1(z+(j_0-2)c)\cdots Q_1(z+c)\eeas
if $j_0\geq 2$ and 
\beas U(z,Q_1(z))&=&\sum\limits_{2\le j\leq n}c_j\frac{H(z+jc)}{H(z)}Q_1(z+(j-1)c)Q_1(z+(j-2)c)\cdots Q_1(z+c)\\&&+B_{j_0}(z)\eeas
if $j_0=1$.

Note that 
\[\mu(e^{R_1})=\rho(e^{R_1})=\deg(R_1)=s-1\]
and 
\[\rho(e^{Q-R_{j_0}})=\deg(Q-R_{j_0})\leq s-2<s-1=\mu(e^{R_1}).\]

Consequently 
\bea\label{al.29} T(r,e^{Q-R_{j_0}})=S(r,e^{R_1})=S(r,Q_1).\eea
 
Also from (\ref{al.8}), it is easy to prove that
\bea\label{al.30} m\left(r,\frac{H(z+jc)}{H(z)}\right)=S(r,e^{R_1})=S(r,Q_1)\;(j=1,2,\ldots,n).\eea

Now from (\ref{al.27}), (\ref{al.29}) and (\ref{al.30}), we see that
\beas m(r,B_{j_0}(z))\leq m\left(r,\frac{H(z+j_0c)}{H(z)}\right)+m\left(r,e^{Q(z)-R_{j_0}(z)}\right)\leq S(r,Q_1).\eeas

On the other hand, from (\ref{al.28}), we see that $U(z,Q_1)\not\equiv 0$ and $\deg(U(z,Q_1))=n-1\geq 1$.
Now using Lemma \ref{l5} to (\ref{al.28}), we get $m(r,Q_1)=S(r,Q_1)$ and so $T(r,Q_1)=m(r,Q_1)=S(r,Q_1)$, which is impossible.\par

\medskip
{\bf Sub-case 2.2.2.} Let $jsa_sc\neq d_{s-1}$, where $1\leq j\leq n$. Now (\ref{al.24}) can be rewriten as 
\bea\label{al.31} e^{Q(z)}=e^{d_{s-1}z^{s-1}}e^{\tilde  P_{s-2}(z)}=\sum\limits_{j=0}^nc_j\frac{H(z+jc)}{H(z)}e^{R_j(z)},\eea
where 
\bea\label{al.32} \tilde P_{s-2}(z)=Q(z)-d_{s-1}z^{s-1}=d_{s-2}z^{s-2}+d_{s-3}z^{s-3}+\cdots+d_0.\eea 

Again from (\ref{al.24a}) and (\ref{al.31}), we have
\bea\label{al.33} e^{Q(z)}=e^{d_{s-1}z^{s-1}}e^{\tilde  P_{s-2}(z)}=\sum\limits_{j=1}^nc_j\frac{H(z+jc)}{H(z)}e^{jsa_scz^{s-1}} e^{P_{s-2,j}(z)}+(-1)^n.\eea

Note that 
\[ns|a_sc|>(n-1)s|ca_s|>\cdots>s|a_sc|\]
and either $|d_{s-1}|\in\{js|a_sc|: j=1,2,\ldots,n\}$ or $|d_{s-1}|\not\in\{js|a_sc|: j=1,2,\ldots,n\}$.
Therefore if we compare $|d_{s-1}|$ with $ns|a_sc|$, $(n-1)s|a_sc|$, $\cdots$, $s|a_sc|$, then it is enough to compare $|d_{s-1}|$ with $ns|a_sc|$.

Now we distinguish between the following two sub-cases.\par

\medskip
{\bf Sub-case 2.2.2.1.} Let $ns|a_sc|\leq |d_{s-1}|$. Suppose $\arg d_{s-1}=\theta _1$ and $\arg (a_sc)=\theta _2$. Take $\theta_0$ such that $\cos ((s-1)\theta_0+\theta_1)=1$. Then by Lemma \ref{l7}, we see that for any given $\varepsilon\;(0<\varepsilon<s-\rho(H))$, there exists a set $E\subset (1,\infty)$, where $l(E)<+\infty$ such that for all $z=re^{i\theta_0}$ satisfying $|z|=r\not\in [0,1]\cup E$, we have
\bea\label{al.34} \exp\left(-r^{\rho(H)-1+\varepsilon}\right)\leq \bigg|\frac{H(z+jc)}{H(z)}\bigg|\leq \exp\left(r^{\rho(H)-1+\varepsilon}\right)\;(j=1,2,\ldots,n).\eea

Note that 
\bea\label{al.35} &&\left|\exp\left(d_{s-1}z^{s-1}\right)\right|\\&=&\left|\exp\left(|d_{s-1}|r^{s-1}\left(\cos((s-1)\theta_0+\theta_1))+i\sin((s-1)\theta_0+\theta_1)\right)\right)\right|\nonumber=\exp\left(|d_{s-1}|r^{s-1}\right).\nonumber\eea

Similarly we can show that
\bea\label{al.36} \left|\exp\left(jsa_scz^{s-1}\right)\right|=\exp\left(js|a_{s}c|r^{s-1}\cos((s-1)\theta_0+\theta_2)\right),\;j=1,2,\ldots,n.\eea

If possible suppose $\deg(P_{s-2,j})=s-2>0$. Set $P_{s-2,j}=c_{0,j}+c_{1,j}z+\ldots+c_{s-2,j}z^{s-2}$,
where $c_{k,j}=\varrho_{k,j}(\cos \alpha_{k,j}+i\sin\alpha_{k,j})$. Then
\beas P_{s-2,j}(z)=\sum\limits_{k=0}^{s-2}\varrho_{k,j}r^{k}\{\cos(\alpha_{k,j}+k\theta_0)+i\sin(\alpha_{k,j}+k\theta_0)\}\eeas
and so
\bea\label{ax.1} \left|\exp\left(P_{s-2,j}(z)\right)\right|=\exp\left(\sum\limits_{k=0}^{s-2}\varrho_{k,j}r^{k}\cos(\alpha_{k,j}+k\theta_0)\right).\eea

Note that 
\beas \left|\sum\limits_{k=0}^{s-2}\varrho_{k,j}r^{k}\cos(\alpha_{k,j}+k\theta_0)\right|\leq \sum\limits_{k=0}^{s-2}\varrho_{k,j}r^{k}=\varrho_{s-2,j}r^{s-2}\Big(1+\frac{\varrho_{s-3,j}}{\varrho_{s-2,j}}\frac{1}{r}+\ldots+\frac{\varrho_{0,j}}{\varrho_{s-2,j}}\frac{1}{r^{s-2}}\Big).\eeas

We see that $\frac{\varrho_{s-3,j}}{\varrho_{s-2,j}}\frac{1}{r}+\ldots+\frac{\varrho_{0,j}}{\varrho_{s-2,j}}\frac{1}{r^{s-2}}\rightarrow 0$ as $r\rightarrow \infty$
and so we may assume that
\[\frac{\varrho_{s-3,j}}{\varrho_{s-2,j}}\frac{1}{r}+\ldots+\frac{\varrho_{0,j}}{\varrho_{s-2,j}}\frac{1}{r^{s-2}}<\frac{1}{2}\]
for sufficiently large values of $r$. Let 
\[|C_0|=\max\{|c_{s-2,j}|: j=1,\ldots,n\}.\]

Consequently from (\ref{ax.1}), we get
\bea\label{al.37} \left|\exp\left(P_{s-2,j}(z)\right)\right|<\exp\left(\frac{3|C_0|}{2}r^{s-2}\right),\;j=1,2,\ldots,n\eea
for sufficiently large values of $r$. If $P_{s-2,j}$ is a constant, then (\ref{al.37}) also holds.

Similarly, we can prove that 
\bea\label{al.38} \left|\exp\left(-\tilde P_{s-2}(z)\right)\right|<\exp\left(\frac{3|d_{s-2}|}{2}r^{s-2}\right)\eea
for sufficiently large values of $r$.

Clearly from (\ref{al.34}), (\ref{al.36}) and (\ref{al.37}), we get 
\bea\label{al.39} &&\left|\frac{H(z+jc)}{H(z)}e^{jsa_scz^{s-1}} e^{P_{s-2,j}(z)}\right|\\&\leq& \exp\left(js|a_{s}c|r^{s-1}\cos((s-1)\theta_0+\theta_2))+r^{\rho(H)-1+\varepsilon}+\frac{3|C_0|}{2}r^{s-2}\right)\nonumber\\&=&\exp (jx+y),\nonumber\eea
where 
\[x=s|a_{s}c|r^{s-1}\cos((s-1)\theta_0+\theta_2))\;\;\text{and}\;\;y=r^{\rho(H)-1+\varepsilon}+\frac{3|C_0|}{2}r^{s-2},\] 
where $j=1,2,\ldots,n$.

Now from (\ref{al.33}), (\ref{al.35}), (\ref{al.38}) and (\ref{al.39}), we have
\bea\label{al.40} \exp\left(|d_{s-1}|r^{s-1}\right)&=&\left|\exp(d_{s-1}z^{s-1})\right|
\\&=&\left|\frac{\exp(Q(z))}{\exp(\tilde P_{s-2}(z))}\right|\nonumber\\
&\leq & \left|\sum\limits_{j=1}^nc_j\frac{H(z+jc)}{H(z)}e^{jsa_scz^{s-1}} e^{P_{s-2,j}(z)}+(-1)^n\right|\left|\exp(-\tilde P_{s-2}(z))\right|\nonumber\\&\leq &
\left(\sum\limits_{j=1}^n |c_j|\exp \left(jx+y\right)+1\right)\exp\left(\frac{3|d_{s-2}|}{2}r^{s-2}\right)
\nonumber.\eea


\medskip
First we suppose that $x\leq 0$, i.e., $\cos((s-1)\theta_0+\theta_2))\leq 0$. Let $\cos((s-1)\theta_0+\theta_2))=-x_1$, where $0\leq x_1\leq 1$. If $x_1=0$, then from (\ref{al.40}), we get
\bea\label{ax.2} \exp\left(|d_{s-1}|r^{s-1}\right)<C_1\exp\left(r^{\rho(H)-1+\varepsilon}+\frac{3(|C_0|+|d_{s-2}|)}{2}r^{s-2}\right),\eea
where $C_1>0$. Since $\rho(H)-1+\varepsilon<s-1$, (\ref{ax.2}) leads to a contradiction. Hence $0<x_1\leq 1$.
In this case, since $\rho(H)-1+\varepsilon<s-1$, we see that
\[\sum\limits_{j=1}^n |c_j|\exp \left(-js|a_{s}c|x_1r^{s-1}+r^{\rho(H)-1+\varepsilon}+\frac{3|C_0|}{2}r^{s-2}\right)\rightarrow 0\] as $r\rightarrow \infty$
and so for $\varepsilon_1$ in $(0,1)$ and $r>r(\varepsilon_1)$, we have 
\[\sum\limits_{j=1}^n |c_j|\exp \left(-js|a_{s}c|x_1r^{s-1}+r^{\rho(H)-1+\varepsilon}+\frac{3|C_0|}{2}r^{s-2}\right)<\varepsilon_1.\]
Consequently from (\ref{al.40}), we get 
\beas \exp\left(|d_{s-1}|r^{s-1}\right)<(1+\varepsilon_1)\exp\left(\frac{3|d_{s-2}|}{2}r^{s-2}\right),\eeas
which leads to a contradiction.

\medskip
Next we suppose that $x>0$, i.e., $\cos((s-1)\theta_0+\theta_2))>0$. In this case, we have $\exp (jx+y)<\exp (nx+y)$
for $j=0,\ldots,n$ and so from (\ref{al.40}), we get
\bea\label{ax.4} &&\exp\left(|d_{s-1}|r^{s-1}\right)\\&\leq &
\left(\sum\limits_{j=1}^n |c_j|\exp \left(jx+y\right)+1\right)\exp\left(\frac{3|d_{s-2}|}{2}r^{s-2}\right)
\nonumber\\&\leq&
\left(\sum\limits_{j=0}^n|c_j|\right)\exp \left(jx+y\right)\exp\left(\frac{3|d_{s-2}|}{2}r^{s-2}\right)
\nonumber\\
&=&B\exp\left(ns|a_{s}c|r^{s-1}\cos((s-1)\theta_0+\theta_2))+r^{\rho(H)-1+\varepsilon}+\frac{3(|C_0|+|d_{s-2}|)}{2}r^{s-2}\right)\nonumber,\eea
where 
\[B=\sum\limits_{j=0}^n|c_j|.\]

Since 
\[\rho(H)-1+\varepsilon<s-1\;\;\text{and}\;\;B=\exp(\log B)=o(r^{s-1}),\]
from (\ref{ax.4}), we deduce that
\bea\label{al.41} \exp\left(|d_{s-1}|r^{s-1}\right)\leq \exp\left(ns|a_sc|\cos((s-1)\theta_0+\theta_2)r^{s-1}+o(r^{s-1})\right).\eea

By assumption, we have $d_{s-1}\neq nsa_sc$ and $ns|a_sc|\leq |d_{s-1}|$. 

If $ ns|a_sc|=|d_{s-1}|$, then $\cos((s-1)\theta _0+\theta _2)\neq 1$ and so $\cos((s-1)\theta _0+\theta _2)<1$. Therefore $ns|a_sc|\cos((s-1)\theta _0+\theta _2)<ns|a_sc|=|d_{s-1}|$.

If $ns|a_sc|<|d_{s-1}|$, then obviously $ns|a_sc|\cos((s-1)\theta _0+\theta _2)\leq ns|a_sc|<|d_{s-1}|$. 

\smallskip
Therefore in either case, we have 
\[ns|a_sc|\cos((s-1)\theta _0+\theta _2)<|d_{s-1}|.\]

Then there exists $\varepsilon_2>0$ such that $ns|a_sc|\cos((s-1)\theta _0+\theta _2)+2\varepsilon_2<|d_{s-1}|$ and so from (\ref{al.41}), we have 
\beas \exp\left(|d_{s-1}|r^{s-1}\right)\leq \exp\left(\left(|d_{s-1}|-2\varepsilon_2\right)r^{s-1}\right),\eeas
which is a contradiction.\par

\medskip
{\bf Sub-case 2.2.2.2.} Let $ns|a_sc|>|d_{s-1}|$. In this case, (\ref{al.33}) can be rewritten as 
\bea\label{xa.1} &&e^{nsa_scz^{s-1}} e^{P_{s-2,n}(z)}\\&=&\sum\limits_{j=0}^{n-1}\frac{c_j}{c_n}\frac{H(z+jc)}{H(z+nc)}e^{jsa_scz^{s-1}} e^{P_{s-2,j}(z)}-
\frac{1}{c_n}\frac{H(z)}{H(z+nc)}e^{d_{s-1}z^{s-1}}e^{\tilde  P_{s-2}(z)}\nonumber.\eea

Let $\arg (a_sc)=\theta _1$ and $\arg d_{s-1}=\theta _2$. Take $\theta_0$ such that $\cos ((s-1)\theta_0+\theta_1)=1$. Then using
Lemma \ref{l7}, we see that for any given $\varepsilon\;(0<\varepsilon<s-\rho(H))$, there exists a set $E\subset (1,\infty)$, where $l(E)<+\infty$ such that for all $z=re^{i\theta_0}$ satisfying $|z|=r\not\in [0,1]\cup E$, we have
\bea\label{xa.2} \exp\left(-r^{\rho(H)-1+\varepsilon}\right)\leq \bigg|\frac{H(z+jc)}{H(z+nc)}\bigg|\leq \exp\left(r^{\rho(H)-1+\varepsilon}\right)\;(j=0,1,\ldots,n-1).\eea

We easily get
\bea\label{xa.3} \left|\exp\left(jsa_scz^{s-1}\right)\right|=\exp\left(js|a_{s}c|r^{s-1}\right),\;j=0,1,\ldots,n\eea
and
\bea\label{xa.4} \left|\exp\left(d_{s-1}z^{s-1}\right)\right|=\left|\exp\left(|d_{s-1}|r^{s-1}\cos((s-1)\theta_0+\theta_2)\right)\right|\leq \left|\exp\left(|d_{s-1}|r^{s-1}\right)\right|.\eea

Also we easily obtain
\bea\label{xa.5} \left|\exp (\pm P_{s-2,j}(z))\right|<\exp\Big(\frac{3|C_0|}{2}r^{s-2}\Big)\;\text{and}\;\left|\exp (\tilde P_{s-2}(z))\right|<\exp\Big(\frac{3|d_{s-2}|}{2}r^{s-2}\Big)\eea
for sufficiently large values of $r$. 

Then from (\ref{xa.2}), (\ref{xa.3}) and (\ref{xa.5}), we get
\bea\label{xa.6} \left|\frac{H(z+jc)}{H(z+nc)}e^{jsa_scz^{s-1}} e^{P_{s-2,j}(z)}\right| &\leq& \exp\left(js|a_{s}c|r^{s-1}+r^{\rho(H)-1+\varepsilon}+\frac{3|C_0|}{2}r^{s-2}\right)\\&\leq & \exp\left((n-1)s|a_{s}c|r^{s-1}+r^{\rho(H)-1+\varepsilon}+\frac{3|C_0|}{2}r^{s-2}\right),\nonumber\eea
for $j=1,\ldots, n-1$. 

Again from (\ref{xa.2}), (\ref{xa.4}) and (\ref{xa.5}), we get
\bea\label{xa.7} \left|\frac{H(z)}{H(z+nc)}e^{d_{s-1}z^{s-1}}e^{\tilde  P_{s-2}(z)}\right|\leq \exp\left(|d_{s-1}|r^{s-1}+r^{\rho(H)-1+\varepsilon}+\frac{3|d_{s-2}|}{2}r^{s-2}\right).\eea

Now from (\ref{xa.1}), (\ref{xa.3}), (\ref{xa.5})-(\ref{xa.7}), we obtain
\bea\label{xa.8}  &&\exp\left(ns|a_{s}c|r^{s-1}\right)\\&=&\left|\sum\limits_{j=0}^{n-1}\frac{c_j}{c_n}\frac{H(z+jc)}{H(z+nc)}e^{jsa_scz^{s-1}} e^{P_{s-2,j}(z)}-
\frac{1}{c_n}\frac{H(z)}{H(z+nc)}e^{d_{s-1}z^{s-1}}e^{\tilde  P_{s-2}(z)}\right|\left|e^{P_{s-2,n}(z)}\right|\nonumber\\&\leq&
B\exp\left(x_2r^{s-1}+r^{\rho(H)-1+\varepsilon}+\frac{6|C_0|+3|d_{s-2}|}{2}r^{s-2}\right),\nonumber
\eea
where 
\[B=\sum\limits_{j=0}^n |c_j|\;\;\text{and}\;\;x_2=\max\{(n-1)s|a_{s}c|, |d_{s-1}|\}.\]

Since $|d_{s-1}|<ns|a_{s}c|$, we have $x_2<ns|a_{s}c|$. Also we have $\rho(H)-1+\varepsilon<s-1$. Consequently from (\ref{xa.8}), we get a contradiction.

This completes the proof. 
\end{proof}

\begin{proof}[{\bf Proof of Theorem \ref{t2}}]
By the given conditions, we have $\lambda(f)<\rho(f)$. Then there exist an entire function $H(\not\equiv 0)$ and a polynomial $P$ such that $f=He^{P}$, where $\rho(H)<\rho(f)$ and $\deg(P)=\rho(f)$. 

Now we divide the following two cases.\par

\medskip
{\bf Case 1.} Let $\rho(f)<2$. Since $\rho(f)=\deg(P)<2$, it follows that $\rho(f)=\deg(P)=1$. Again since $\lambda(H)<\rho(f)$, we have $\lambda(H)<1$. Therefore we may assume that 
\bea\label{all.2} f(z)=H(z)e^{dz},\eea
where $d\in\mathbb{C}\backslash \{0\}$ and $\lambda(H)<1$. Now from (\ref{all.2}), we get
\bea\label{all.3} \Delta_c^nf(z)=H_n(z)e^{dz},\eea
where 
\bea\label{all.4} H_n(z)=\sum\limits_{i=0}^nc_iH(z+ic)e^{icd}\;\;\text{and}\;\;c_i=(-1)^{n-i}C^{i}_{n}.\eea

Since $f$ and $\Delta_c^nf$ share $0$ CM, then by \emph{Theorem B}, we have  
\bea\label{all.4}\Delta_c^nf=af,\eea
where $a\in\mathbb{C}\backslash \{0\}$. Then from (\ref{all.2})-(\ref{all.4}), we get 
\bea\label{all.5} \sum\limits_{i=1}^nc_{i}e^{icd}H(z+ic)+((-1)^n-a)H(z)=0.\eea

Now we divide the following two sub-cases.\par

\medskip
{\bf Sub-case 1.1.} Let $a\neq (-1)^n$. Now using Lemma \ref{l8} to (\ref{all.5}), we conclude that
$H$ is a non-zero polynomial. We claim that $H$ is a non-zero constant. If not, suppose $\deg(H)=k\in\mathbb{N}$. Let $H(z)=a_kz^k+a_{k-1}z^{k-1}+\cdots+a_1z+a_0\;(a_k\neq 0)$. Then comparing the coefficients of $z^k$ and $z^{k-1}$ from (\ref{all.5}), we have respectively
\bea\label{all.6}\sum\limits_{i=1}^n c_{i}e^{icd}+((-1)^n-a)=(e^{cd}-1)^n-a=0\eea
and
\bea\label{all.7}\sum\limits_{i=1}^n c_{i}e^{icd}(ika_k+a_{k-1})+((-1)^n-a)a_{k-1}=0.\eea 
 
Now from (\ref{all.6}) and (\ref{all.7}), we immediately have
\bea\label{all.8}\sum\limits_{i=1}^n c_{i}ie^{icd}=0.\eea

Again from (\ref{all.6}), we have $a=(e^{cd}-1)^n$. Since $a\neq 0$, we have $e^{cd}\neq 1$.
Note that 
\beas C_n^ii=\frac{n!i}{(n-i)!i!}=\frac{n(n-1)!}{((n-1)-(i-1))!(i-1)!}=nC_{n-1}^{i-1}\eeas
and so
\bea\label{all.9}\sum\limits_{i=1}^n c_{i}ie^{icd}=\sum\limits_{i=1}^n (-1)^{n-i}C^{i}_{n}ie^{icd}&=&ne^{cd}\sum\limits_{i=1}^n(-1)^{n-i}C^{i-1}_{n-1}e^{(i-1)cd}\\&=&ne^{cd}(e^{cd}-1)^{n-1}\neq 0.\nonumber\eea

Therefore from (\ref{all.8}) and (\ref{all.9}), we get a contradiction.
Hence $H$ is a non-zero constant. Then from (\ref{all.2}), we may assume that $f(z)=c_0e^{dz}$, where $c_0, d\in\mathbb{C}\backslash \{0\}$.\par

\medskip
{\bf Sub-case 1.2.} Let $a=(-1)^n$. Then from (\ref{all.5}), we have 
\bea\label{all.10}\sum\limits_{i=1}^n c_{i}e^{icd}H(z+ic)=0.\eea

If $n=1$ then we immediately get a contradiction from (\ref{all.10}). Hence $n\geq 2$.

If we let $w=z+c$, then from (\ref{all.10}), we have
\bea\label{all.11} \sum\limits_{i=0}^{n-1} c_{i+1}e^{(i+1)cd}H(w+ic)=0.\eea

Note that $\rho(H)<1$. Now using Lemma \ref{l8} to (\ref{all.11}), we conclude that $H(w)$ is a non-zero polynomial, i.e., $H(z+c)$ is non-zero polynomial. Consequently, $H(z)$ is a non-zero polynomial. We claim that $H$ is a non-zero constant. If not, suppose $\deg(H)=m\in\mathbb{N}$. Let $H=b_mz^m+b_{m-1}z^{m-1}+\cdots+b_0(b_m\neq 0)$. 
Now comparing the coefficient of $z^m$ from (\ref{all.11}), we have 
 \bea\label{all.12}\sum\limits_{i=1}^nc_ie^{icd}=0.\eea
 
Again comparing the coefficient of $z^{m-1}$ from (\ref{all.11}) and then using (\ref{all.12}), we have
\[\sum\limits_{i=1}^nc_ie^{icd}(mib_m+b_{m-1})=0,\]
i.e.,
\bea\label{all.13}mb_m\sum\limits_{i=1}^nc_iie^{icd}=0.\eea

On the other hand form (\ref{all.12}), we have
\[\sum\limits_{i=1}^n(-1)^{n-i}C^{i}_{n}e^{icd}+(-1)^n-a=0,\]
i.e.,
\bea\label{all.14}a=(e^{cd}-1)^n\neq 0.\eea

Consequently from (\ref{all.9}), (\ref{all.13}) and  (\ref{all.14}), we get a contradiction. So $H$ is a non-zero constant.
Hence $f$ is of the form $f(z)=c_0e^{dz}$, where $c_0, d\in\mathbb{C}\backslash \{0\}$.\par

\medskip
{\bf Case 2.} Let $\rho(f)\geq 2$. In this case, proceeding in the same way as done in \emph{Case 2} in the proof of Theorem \ref{t1}, we get a contradiction. 

Hence the proof is complete.
\end{proof}

\section{\bf {Statements and declarations}}
\vspace{1.3mm}

\noindent \textbf {Conflict of interest:} The authors declare that there are no conflicts of interest regarding the publication of this paper.\vspace{1.5mm}

\noindent{\bf Funding:} There is no funding received from any organizations for this research work.
\vspace{1.5mm}

\noindent \textbf {Data availability statement:}  Data sharing is not applicable to this article as no database were generated or analyzed during the current study.

\end{document}